\documentclass[twoside,12pt,a4paper,notitlepage]{amsart}
\usepackage{amssymb}

\begin{document}

\title[The Quasi-Stationary Distribution of a SIS model]
 {Approximations of the Quasi-Stationary Distribution of a Logistic SIS Model for Endemic Infections}

\author{Ingemar N{\aa}sell}
\address{Department of Mathematics \\
        The Royal Institute of Technology \\
        S-100 44 Stockholm, Sweden}
\email{ingemar@kth.se}

 
\date{\today}

\begin{abstract}
Errors of approximations of the quasi-stationary distribution (the QSD) of the logistic SIS model are evaluated numerically. 
The results are used to derive asymptotic approximations of the approximation errors for large populations. 
We show in particular that there are two approximations above threshold for which the approximation errors are exponentially small. 
One of these approximations has been known for some time, while the other one is new. 
The result that the older one of these two approximations has an exponentially small approximation error is new. 
\end{abstract}

\maketitle

\section{Introduction}

The stochastic SIS model that we study is a continuous time Markov Chain with a finite state space. 
It serves to model the number of infected individuals $I(t)$ in a constant population of $N$ hosts, where individuals that recover from infection are immediately susceptible to new infections.  
This model was first dealt with by Weiss and Dishon (1971).  
As described by Kryscio and Lefèvre (1989), the same model has been used for other applications, such as the propagation of rumours, a particular chemical reaction, and the growth of a population with limiting growth factors. 

The behaviour of the model is determined by the birth rate $\lambda_n$ and the death rate $\mu_n$, defined as follows:  
\begin{align} \label{1.1}
    & \lambda_n = \mu R_0 \left(1 - \frac{n}{N}\right) n, \quad n=0,1,\dots,N, \\ \label{1.2}
    & \mu_n = \mu n, \quad  n=0,1,\dots,N.
\end{align}

The parameter space contans three parameters, namely the poulation size $N$, the threshold parameter $R_0$, and the recovery rate $\mu$. 
Among these parameters, $N$ is a dimensionless positive integer, $R_0$ is a dimensionless 
positive threshold parameter, and $\mu$ is a positive number with the dimension inverse time. 

The model that we study is a special case of the Verhulst logistic model. 
The deterministic version of this classical model was formulated by Verhulst (1838). 
The stochastic version of this model is a finite-state birth-death process with the same 
birth-rate as in \eqref{1.1}, while its death-rate is written 
\begin{equation} \label{1.3}
    \mu_n = \mu \left(1 + \alpha \frac{n}{N} \right)n, \quad n=0,1,\dots,N,
\end{equation}
with $\alpha\ge 0$. 
Thus, the SIS model is a Verhulst logistic model with $\alpha=0$. 

It is important to recognize that the SIS model has an absorbing state at the origin.  
We study the number of infected individuals in the case when absorption has not taken place. 
This is done by conditioning the state variable $I(t)$ on non-extinction. 
The stationary distribution of the conditioned state variable $I^{(Q)}(t)$ is the so-called quasi-stationary distribution (QSD). 
It has been studied extensively, as shown e.g. by Cavender (1978), Kryscio and Lefèvre (1989),  
N{\aa}sell (1996), (1999), (2011), Ovaskainen (2001), Clancy and Mendy (2010), and Clancy (2012). 
Exact expressions for the QSD are not available. 
This leads to a search for useful approximations.

It is straightforward to show that the model just described has a threshold at $R_0=1$. 
The deterministic version of the model predicts that any infection present at $t=0$ will ultimately approach a positive steady-state infection level for $R_0>1$, while it will ultimately disappear if 
$R_0 \le 1$. 
These are clear qualitative differrences. 
In accordance with this, we expect the stochastic version of the model to also show qualitatively different behaviours above and below threshold. 
Clearly, this means that the QSD, and also approximations of the QSD, are entirely different above and below threshold. 
Kryscio and Lefèvre (1989) give two explicit distributions, where one serves to approximate the QSD above threshold, and the other one below threshold.  
We describe these two approximations in Section 2.

A different approximation of the QSD above threshold is given in Section 3. 
It is based on an approximation derived by Ovaskainen (2001). 
The derivation of this approximation actually takes the form of a three-stage modification of the 
Ovaskainen result.

Several approximations of the QSD of the SIS model in addition to the three that are described in Sections 2 and 3 have been suggested. 
Examples are given by Clancy and Mendy (2010).  
We compare competing approximations of the QSD via a determination of their approximation errors. 
We study the error for any approximation of the QSD as a function of the population size $N$ with fixed $R_0$, and search for the magnitude of the approximation error as a function of $N$. 
Numerial evaluations are then necessary, since no explicit expression is available for the QSD. 

Section 4 is used to report the results of numerical evaluations of approximation errors for several 
approximations of the QSD. 
We give results both above and below the threshold. 
We note that the quasi-stationary distribution is of possibly less interest below the threshold, since then the time to extinction is so short that there may not be enough time for the quasi-stationary distribution to establish itself. 

The paper ends with some concluding comments in Section 5.

\section{Stationary Distributions of Two Auxiliary Processes} 

The two approximations of the QSD derived by Kryscio and Lefèvre (1989) are derived here with a slightly different approach. 
Our notation is also different; we use the same notation as in N{\aa}sell (2011).

We introduce two auxiliary processes $I^{(0)}(t)$ and $I^{(1)}(t)$ to describe our approach.
Both of them are birth-death processes that are similar to the SIS model that we study, but with the important differences that they lack absorbing states.  
They were introduced by Cavender (1978) and Kryscio and Lefevre (1989), respectively. 
The state spaces of both of them are equal to the state space $\{1,2,\dots,N\}$ of the conditioned state variable $I^{(Q)}(t)$. 
The birth rates of both of the auxiliary processes are equal to the birth rates of the SIS process, while the death rates differ slightly. 
Thus, the death rates of the $I^{(0)}(t)$-process are equal to the death rates of the SIS process, with the one important exception that the death rate $\mu_1$ is replaced by zero. 
Furthermore, the death rates of the $I^{(1)}(t)$-process are all smaller than the death rates of the 
SIS process, with $\mu_n$ replaced by $\mu_{n-1}$.

The stationary distributions of the two auxiliary processes are denoted by
$p^{(0)} =\left(p_1^{(0)},p_2^{(0)},\dots,p_N^{(0)}\right)$ and    
$p^{(1)} =\left(p_1^{(1)},p_2^{(1)},\dots,p_N^{(1)}\right)$, respectively. 
We determine explicit expressions for them.
Using results in N{\aa}sell (2011), we find that they are equal to  
\begin{align} \label{2.1}
   & p_n^{(0)} = \frac{\pi_n}{\sum_{n=1}^N \pi_n}, \quad n=1,2,\dots,N, \\ \label{2.2}
   & p_n^{(1)} = \frac{\rho_n}{\sum_{n=1}^N \rho_n}, \quad n=1,2,\dots,N,
\end{align}
where 
\begin{align} \label{2.3}
    & \pi_1=1, \quad \pi_n = \frac{\lambda_1 \lambda_2 \cdots \lambda_{n-1}}
{\mu_2 \mu_3 \cdots \mu_n}, \quad n=2,3,\dots,N, \\ \label{2.4}
    & \rho_1 = 1, \quad \rho_n = \frac{\lambda_1 \lambda_2 \cdots \lambda_{n-1}}
{\mu_1 \mu_2 \cdots  \mu_{n-1}}, \quad n =2,3,\dots,N.
\end{align}
By using the expressions for $\lambda_n$ and $\mu_n$ in \eqref{1.1} and \eqref{1.2} we find that 
$\pi_n$ and $\rho_n$ can be written as follows: 
\begin{align} \label{2.5}
    & \pi_n = \frac{1}{n} \frac{1}{R_0} \frac{N!}{(N-n)!} \left(\frac{R_0}{N} \right)^n, 
         \quad n=1,2,\dots,N, \\ \label{2.6}
   & \rho_n =  \frac{1}{R_0} \frac{N!}{(N-n)!} \left(\frac{R_0}{N}\right)^n, 
         \quad n=1,2,\dots,N. 
\end{align}

An explicit expression for the stationary ditribution $p^{(0)}$ can now be found by inserting the expression \eqref{2.5} for $\pi_n$ into \eqref{2.1}.
In similarity to this we determine an explicit expression for the stationary distribution $p^{(1)}$ 
by inserting the expression \eqref{2.6} for $\rho_n$ into \eqref{2.2}. 
Both of the two stationary distributions $p^{(0)}$ and $p^{(1)}$ can be used to approximate the QSD, both above and below threshold. 
The magnitudes of the resulting approximation errors are determined using numerical evaluations in Section 4. 
We show there that $p^{(0)}$ is the preferred approximation above threshold, while  $p^{(1)}$ is preferred below threshold. 
The facts that the QSD is approximated by  $p^{(0)}$ above threshold and by  $p^{(1)}$ below threshold are not new. 
They were given by Kryscio and Lefèvre (1989).


\section{Ovaskainen's Approximation, and Modifications}

Ovaskainen (2001) presents an approximation of the QSD of the SIS model that we deal with. 
In this work, he uses a different parametrization than above. 
He uses one parameter that he calls $R_0$. 
Since it is different from our parameter with the same name, we rename his $R_0$ and call it $R_0^{(OV)}$. 
It then turns out that his parameter $R_0^{(OV)}$ can be expressed as follows in terms of our parameters $N$ and $R_0$:
\begin{equation} \label{3.1}
    R_0^{(OV)} = \frac{N-1}{N} R_0.
\end{equation}
We proceed to describe Ovaskainen's result, and three modifications of it.

Ovaskainen gives two theorems that both are claimed to give approximations of the QSD under certain conditions. 
We denote these approxmations by $q^{(OVa)}$ and $q^{(OVb)}$, respectively. 
To express them we introduce $F_i(N,R_0)$ as  follows:
\begin{equation} \label{3.2}
    F_i(N,R_0) =  \left (1 - \frac{N}{(N-1)R_0}\right)\exp\left(-\frac{N}{R_0}\right) 
        \left[1 - \left(\frac{N}{(N-1)R_0}\right)^i  \right].
\end{equation}
The two approximations can then be written as follows:
\begin{equation} \label{3.3} 
     q_i^{(OVa)} = F_i(N,R_0) q_i^{(OVb)}, 
        \quad i=1,2,\dots,N, \quad R_0>1, \quad N\to\infty,
\end{equation}
and 
\begin{multline} \label{3.4}
    q_i^{(OVb)} = \frac{N}{i(N-i)!}\left(\frac{N}{R_0}\right)^{N-i}, \quad i=1,2,\dots,N-1, \\
     \quad N\ge 2, \quad R_0\to\infty,
\end{multline}
\begin{equation} \label{3.5}
   q_N^{(OVb)} = 1 - \frac{N^2}{(N-1)R_0},  \quad N\ge 2, \quad R_0\to\infty,
\end{equation}

Ovaskainen has used numerical experimentation to show that adjustments in these approximations are needed when $N$ is small. 
The need for adjustments arises when the function $K$ defined as follows is larger than 1:
\begin{equation} \label{3.6}
    K(N,R_0) = \frac{2(N-1)^2 R_0^2}{N[(N-1)R_0-N]^2}
\end{equation}
In the present paper we consider only cases where $R_0>1$ and $N\to\infty$.
This implies that $K(N,R_0)<1$. 
A consequence of this is that the adjustment required when $K(N,R_0)>1$ plays no role in this paper.

Insertions of the expressions \eqref{3.4} and \eqref{3.5} for $q_i^{(OVb)}$ into \eqref{3.3} 
gives an expression for $q_i^{(OVa)}$ that requires two contradictory conditions for its validity, namely both that $R_0>1$ and $N\to\infty$ and also that $N\ge2$ is fixed and $R_0\to\infty$. 
It is clearly impossible to satisfy these conditions. 
To make progress at this point we introduce a modification of this result. 
It is found by  replacing the conditions that $R_0\to\infty$ and $N \ge 2$ in \eqref{3.3} and \eqref{3.4} by $R_0>1$ and  $N\to\infty$. 
We admit that it is unjustified to introduce this modification. 
We claim however that two additional modifications of the Ovaskainen approximation will lead to interesting results. 
It appears also that this one-stage modification of the Ovaskainen approximation has been used in the numerical results given by Clancy and  Mendy (2010). 
We denote the resulting one-stage modification of the Ovaskainen approximation of the QSD by 
$q^{(OV1)}$.  

It can be expressed as follows:
\begin{multline} \label{3.7}
     q_i^{(OV1)} =  N F_i(N,R_0) \frac{1}{i(N-i)!} \left(\frac{N}{R_0} \right)^{N-i}, \\ 
       i=1,2,\dots,N-1, \quad R_0>1, \quad N \to\infty,
\end{multline}
and
\begin{equation} \label{3.8} 
    q_N^{(OV1)} = F_N(N,R_0) \left(1 - \frac{N^2}{(N-1) R_0} \right), \quad R_0>1, \quad N\to\infty.
\end{equation}

This approximation has the weakness that the probability $q_N^{(OV1)}$ is negative with a large absolute value. 
To avoid this, we introduce a second modification of the Ovaskainen approximation by using the expression \eqref{3.7} also for $i=N$. 

The resulting approximation of the QSD is denoted $q^{(OV2)}$.
It is written as follows:
\begin{multline} \label{3.9}
     q_i^{(OV2)} =  N F_i(N,R_0) \frac{1}{i(N-i)!} \left(\frac{N}{R_0} \right)^{N-i}, \\ 
       i=1,2,\dots,N, \quad R_0>1, \quad N \to\infty.
\end{multline}

To assure that our approximation of the QSD is a true probability distribution, we require the sum of 
the expressions $q_i^{(OV2)}$ over $i$ from 1 to $N$ be equal to 1.  
This leads us to a third approximation step, which is taken by dividing the above expression for $q_i^{(OV2)}$ by the sum of these expressions over $i$ from 1 to $N$.
The resulting expression is denoted $q_i^{(OV3)}$.
Thus, we have 
\begin{equation} \label{3.10} 
     q_i^{(OV3)} = \frac{q_i^{(OV2)}}{\sum_{i=1}^N q_i^{(OV2)}}, 
       \quad i=1,2,\dots,N, \quad R_0>1, \quad N\to\infty. 
\end{equation} 

It is interesting to note that this approximation of the QSD is closely related to the stationary distribution $p^{(0)}$ of the auxiliary process $I^{(0)}(t)$, given by inserting the expression 
\eqref{2.5} for $\pi_n $ into \eqref{2.1}. 
It is straightforward to show that 
\begin{equation} \label{3.11} 
    q_i^{(OV2)} = \pi_i   \left[1-\left(\frac{N}{(N-1)R_0} \right)^i \right], \quad i=1,2,\dots,N.
\end{equation}

\section{Approximation Errors}

Several different distributions have been suggested in the literature as approximations of the QSD 
for the SIS model that we study here. 
It is therefore of interest to be able to compare such approximating distributions.  
As in Clancy and Mendy (2010) we base such comparisons on approximation errors. 
Two different definitions of approximation errors are of interest. 
One, called $Err_1$, is given here,  while a second one, called $Err_2$, was given by Clancy and Mendy (2010). 
Both definitions deal with a QSD written as $q=(q_1,q_2,\dots,q_N)$, and an approximating distribution $\hat{q}=(\hat{q}_1,\hat{q}_2,\dots,\hat{q}_N)$. 
The first error function is defined by 
\begin{equation} \label{4.1}
   Err_1(\hat{q}) = \max_{1\le i \le N}|\hat{q}_i - q_i|.
\end{equation}
and the second one by 
\begin{equation} \label{4.2}
     Err_2(\hat{q}) =  0.5 \sum_{i=1}^N |\hat{q}_i -q_i|,
\end{equation} 
The first error function has a minor advantage over the second one in the sense that the following  inequalities hold: 
\begin{equation} \label{4.3}
   |\hat{q}_i - q_i| \le Err_1(\hat{q}), \quad i=1,2,\dots,N,
\end{equation}
while there are examples that show that the corresponding inequalities do not hold for all 
approximations and all $i$-values if one uses the error function $Err_2$. 
However, both error functions are useful for evaluating magnitudes of errors for approximations 
of the QSD. 
We report below the results of numerical evaluations of approximation errrors $Err_1$ for several approximations of the QSD both above and below the threshold. 

Clancy and Mendy (2010) show values of the approximation error $Err_2$ for several 
approximations of the QSD with several  values of $R_0$ and one constant value of $N$, namely 
$N=50$. 
They include normal, lognormal, binomial, negative binomial, and beta-binomial distributions, and also the one-stage modification of the Ovaskainen result in their study. 
They conclude that among these approximations of the QSD, the beta-binomial distribution is  preferred in a  parameter region that they describe as $R_0>>1$.

We have determined numerical values of the error function $Err_1$ when $R_0>1$ for four different approximations of the QSD, namely the beta-binomial distribution $p^{(B)}$, the stationary distributions $p^{(1)}$ and $p^{(0)}$ of the two auxiliary processes $I^{(1)}(t)$ and $I^{(0)}(t)$, and the three-stage modification $q^{(OV3)}$ of the Ovaskainen  result.
The evaluations have been done for three values of $R_0$ above threshold, namely 2, 5, and 10, and also for three values of $N$, namely 25, 50, and 100.  
It is important for our results to study the approximation error for several $N$-values for each value of $R_0$. 
This will allow us to determine the magnitude of the approximation error.   
The results of our numerical evaluations for $R_0>1$ are listed in Table 1.
All our numerical evaluations have been done using Maple.

\begin{table}[h]
  \begin{center}
    \begin{tabular}{  | c | c | c | c | c | c |  }
       \hline 
  $R_0$ & $N$ & $Err_1(p^{(B)})$ & $Err_1(p^{(1)})$  & $Err_1(p^{(0)})$  & $Err_1(q^{(OV3)})$  \\ \hline
       2   &  25  &  $11*10^{-3}$    &  $23*10^{-3}$  & $7.5*10^{-3}$   & $8.1*10^{-4}$       \\
       2   & 50   &   $4.7*10^{-3}$  &  $11*10^{-3}$   &  $9.0*10^{-5}$ &  $1.5*10^{-5}$  \\ 
       2   & 100 &   $2.1*10^{-3}$  &  $5.2*10^{-3}$  &  $8.1*10^{-9}$  &  $1.4*10^{-9}$   \\ \hline
       5   &  25  &   $24*10^{-4}$  &  $14*10^{-3}$    & $2.7*10^{-9}$  & $1.4*10^{-9}$ \\
       5   &  50  &   $11*10^{-4}$  &  $6.8*10^{-3}$   & $6.1*10^{-18}$ & $2.9*10^{-18}$\\  
       5   & 100  &   $5.0*10^{-4}$ &  $3.3*10^{-3}$  & $2.3*10^{-35}$ & $9.4*10^{-36}$ \\  \hline
     10   &  25   &   $12*10^{-4}$  & $14*10^{-3}$   &  $1.7*10^{-15}$ & $9.8*10^{-16}$ \\
     10   &  50   &   $5.1*10^{-4}$  & $6.3*10^{-3}$  & $1.1*10^{-30}$ & $7.2*10^{-31}$   \\
     10   &  100 &   $2.4*10^{-4}$  & $3.0*10^{-3}$   &$4.6*10^{-61}$  &   $3.0*10^{-61}$\\  \hline
\end{tabular}
   \vskip 4mm
    \caption{Approximation errors for four approximations of the QSD of the SIS model above threshold}
\end{center}
\end{table}

It is seen from Columns 3 and 4 of Table 1 that the two approximation errors that accompany the two approximations $p^{(B)}$ and $p^{(1)}$ are both divided by approximately 2 for each doubling of $N$ when $N$ is sufficiently large. 
This is interpreted as strong indications that these two approximation errors are both of the order 
$O(1/N)$ for large values of $N$. 
We notice also that the beta-binomial distribution $p^{(B)}$ is preferred over the stationary distribution $p^{(1)}$ as an approximation of the QSD when $R_0>1$, since the abserved approximation errors are smaller for the beta-binomial distribution. 
However, considerably smaller approximation errors are reported in Columns 5 and 6 of Table 1. 
It is shown there that the approximation errors that accompany the approximations $p^{(0)}$ and $q^{(OV3)}$ are approximately squared for each doubling of $N$. 
These results are interpreted as strong indications that the corresponding approximation errors are
exponentially small in $N$. 
The two approximations $p^{(0)}$ and $q^{(OV3)}$ for which the approximation errors are given in Columns 5 and 6 of Table 1 are thus found to give considerably smaller approximation errors than the beta-binomial distribution, which was found by Clancy and Mendy to be the preferred approximation for $R_0>>1$. 
A comparison between the magnitudes of the approximation errors in Columns 5 and 6 leads us to conclude that the three-stage modification of the Ovaskainen result is the preferred approximation of the QSD in case $R_0>1$.

\begin{table}[h]
  \begin{center}
    \begin{tabular}{  | c | c | c | c | c | c |  }
       \hline 
  $R_0$ & $N$ & $Err_1\left(p^{(G1)}\right)$ & $Err_1\left(p^{(G2)}\right)$  &
 $Err_1\left(p^{(0)}\right)$  & $Err_1\left(p^{(1)}\right)$  \\ \hline
     0.5 &  25  &  $8.4*10^{-3}$   &  $5.0*10^{-2}$  & $0.19$             & $13*10^{-3}$   \\
     0.5 & 50   &   $4.5*10^{-3}$  &  $2.7*10^{-2}$  &  $0.20$            &  $7.9*10^{-3}$  \\ 
     0.5 & 100 &   $2.4*10^{-3}$  &  $1.4*10^{-2}$  &  $0.21$            &  $4.4*10^{-3}$   \\ \hline
     0.2 &  25  &   $15*10^{-4}$   &  $12*10^{-3}$   & $9.0*10^{-2}$ & $18*10^{-4}$   \\
     0.2 &  50  &   $7.8*10^{-4}$  &  $5.9*10^{-3}$  & $9.3*10^{-2}$ & $9.4*10^{-4}$\\  
     0.2 & 100  &   $4.0*10^{-4}$ &  $3.0*10^{-3}$  & $9.5*10^{-2}$  & $4.8*10^{-4}$ \\  \hline 
     0.1 &  25   &   $3.8*10^{-4}$  & $4.8*10^{-3}$  &  $4.7*10^{-2}$ & $4.1*10^{-4}$ \\
     0.1 &  50   &   $2.0*10^{-4}$  & $2.4*10^{-3}$  & $4.8*10^{-2}$  & $2.1*10^{-4}$   \\
     0.1 &  100 &   $1.0*10^{-4}$  & $1.2*10^{-3}$  & $4.8*10^{-2}$  &  $1.0*10^{-4}$\\  \hline
\end{tabular}
   \vskip 4mm
    \caption{Approximation errors for four approximations of the QSD of the SIS model below threshold} 
\end{center}
\end{table}

Clancy and Mendy (2010) use cumulant closure arguments to derive a geometric distribution as 
an approximation of the QSD for $R_0<1$. 
The first cumulant $\kappa_1$ of this approximation is found to be equal to
\begin{equation} \label{24}
   \kappa_1 = \frac{1}{4} \left( A + \sqrt{A^2 + \frac{8N}{R_0}} \right),
 \end{equation}
where 
\begin{equation} \label{25}
   A = 1 - N \frac{1-R_0}{R_0}.
\end{equation}
They show also that 
\begin{equation} \label{26}
    \kappa_1 = \frac{1}{1-R_0} + \mbox{O}\left(\frac{1}{N}\right) .
\end{equation}
We use  $p^{(G1)}$ denote this geometric distribution. 
We also use $p^{(G2)}$ to refer to the geometric distribution whose expectation equals 
$1/(1-R_0)$. 

We have determined numerical values of the error function $Err_1$ in the parameter region where $R_0<1$ for four different approximations of the QSD. 
The results are shown in Table 2. 
It is seen from Columns 3 and 4 in Table 2 that the approximation errors for the two geometric distributions are both divided by approximately 2 for each doubling of the population size $N$. 
We interpret this as strong indications that these two approximation errors are both of the order 
$O(1/N)$ for large values of $N$. 
We notice also that the geometric distribution that uses the expectation value given by \eqref{24}
gives smaller approximation errors than the slightly simpler one where the expectation is equal to $1/(1-R_0)$. 
The approximation error that accompanies the approximation $p^{(0)}$ of the QSD is seen from 
Column 5 of Table 2 to be practically independent of the population size $N$ when $R_0<1$. 
In contrast to this, we find from Column 6 of Table 2 that the approximation error that goes with the approximation $p^{(1)}$ of the QSD for $R_0<1$ is divided by approximately 2 for each doubling of the population size $N$. 
We interpret this as a strong indication that this approximation error is of the order $O(1/N)$ for large values of $N$. 
A further comparison between the entries in Columns 3 and 6 of Table 2 shows that the approximation errors that accompany the geometric distribution advanced by Clancy and Mendy (2010) are slightly smaller than the approximation errors that are found when the stationary distribution $p^{(1)}$ is used to approximate the QSD for $R_0<1$.

\section{Concluding Comments}

Two important findings in this paper are that both the stationary distribution $p^{(0)}$ of the auxiliary process $I^{(0)}(t)$ and the three-stage modification of the Ovaskainen result provide approximations of the QSD above threshold for which the approximation errors are exponentially small for large values of the population size $N$. 
In addition we note that the latter of these two approximations is superiour to the former, since its approximation errors are smaller than those of the former.  
The results in Table 1 show furthermore that these two approximations of the QSD above threshold 
give approximation errors that are considerably smaller than those that are associated with the beta-binomial distribution, which was found by Clancy and Mendy (2010) to be the best approximant of the QSD above threshold among the several approximations that they consider.   
We summarize by recommending the three-stage modification of the Ovaskainen result as the approximation of the QSD to use above threshold. 
                      
We use the results in Table 2 to recommend the geometric distribution given by Clancy and Mendy (2010) as the approximation of the QSD to use below threshold. 

A comparison between the approximation errors associated with the approximations given by the  stationary distributions $p^{(0)}$ and $p^{(1)}$ of the two auxiliary processes $I^{(0)}(t)$ and $I^{(1)}(t)$ shows that $p^{(0)}$ is the preferred approximation above threshold, while $p^{(1)}$ is the preferred approximation below threshold, in line with results given by Kryscio and Lefèvre (1989). 

We leave it as an open problem to give a rigorous derivation of the approximation of the QSD above threshold that we have described above as a three-stage modification of the Ovaskainen result.

\end{document}